\documentclass[12pt]{article}
\newcommand{\copyleft}{
GNU FDL\thanks{
Copyright (C) 2008 Peter G. Doyle.
Permission is granted to copy, distribute and/or modify this document
under the terms of the GNU Free Documentation License, 
as published by the Free Software Foundation;
with no Invariant Sections, no Front-Cover Texts, and no Back-Cover Texts.
}}
\title{Isospectral hyperbolic surfaces have matching geodesics}
\author{
Peter G. Doyle\thanks{Dartmouth College.}
\and
Juan Pablo Rossetti\thanks{
FaMAF-CIEM, Univ.\ Nac.\ C\'ordoba.
Partially supported by DFG Sonderforschungsbereich 647,
Humboldt University, Berlin.
}}
\date{Version dated 29 April 2008
\\ \copyleft
}
\usepackage{amsrefs}
\newcommand{\Zbl}[1]{Zbl #1}
\newcommand{\half}{\frac{1}{2}}
\newcommand{\wt}{\mathrm{wt}}
\newcommand{\tw}{W}
\newcommand{\trans}{\mathrm{tr}}

\newcommand{\Z}{\mathbf{Z}}
\newcommand{\N}{\mathbf{N}}
\newcommand{\Nonetwo}{\N^+}

\newcommand{\suchthat}{:\,}
\newcommand{\qed}{\rule{2mm}{2.5mm}}
\newcommand{\Podd}{P_\mathrm{odd}}
\newcommand{\intersect}{\cap}
\newcommand{\goesto}{\rightarrow}
\newtheorem{thm}{Theorem}
\newtheorem{prop}{Proposition}
\newtheorem{question}{Question}
\newtheorem{lemma}{Lemma}
\newcommand{\dilemma}[1] {
\left\{
\begin{array}{ll}
#1
\end{array}
\right.
}
\begin{document}

\maketitle

\begin{abstract}
We show that if two closed hyperbolic surfaces (not necessarily orientable 
or even connected) have the same Laplace spectrum,
then for every length they have the same number of
orientation-preserving geodesics
and the same number of orientation-reversing geodesics.
Restricted to orientable surfaces, this result
reduces to Huber's theorem of 1959.
Appropriately generalized,
it extends to hyperbolic 2-orbifolds (possibly
disconnected).
We give examples showing that it fails
for disconnected flat 2-orbifolds.



\end{abstract}

\section{Introduction}

We say that two hyperbolic surfaces
(assumed closed but not necessarily orientable or even connected)
are \emph{almost conjugate}
if their closed geodesics match, in the sense that
for every length $l$ they have the same number of
orientation-preserving geodesics 
and the same number of orientation-reversing geodesics.
Looking beyond dimension 2, we say that two hyperbolic $d$-manifolds
(assumed closed, but not necessarily orientable or even connected)
are almost conjugate
if their geodesics match with respect to length and `twist'.
The twist of a geodesic (also called its `holonomy')
is measured by the conjugacy class in
$O(d-1)$ of the action of parallel translation around the geodesic.
To say that geodesics have matching length and twist amounts to
saying that the corresponding deck transformations are conjugate
under the action of the full isometry group of hyperbolic $d$-space.

This usage of the term `almost conjugate'
accords with Sunada's well-known definition
\cite{sunada:method}
that subgroups
of a finite group are almost conjugate if they meet each conjugacy
class equally.
Here, two connected hyperbolic surfaces are almost conjugate
just if their deck groups intersect
equally every conjugacy class of the full group of isometries of
the hyperbolic plane.
Strictly speaking, it is the deck groups that are almost conjugate,
but it is convenient to apply the term abusively to the surfaces
themselves---the more so since we want to consider disconnected surfaces.

While we haven't specifically required it in the definition,
the matching of geodesics between almost conjugate hyperbolic
manifolds, whether surfaces or manifolds of higher dimension,
can and should be taken to respect the imprimitivity index
of the geodesics as well as their length and twist.

Please note that here and throughout,
by `geodesics' we mean oriented closed geodesics.
Because our geodesics carry a designated orientation,
the number of geodesics of 
length $l$ will always be even, with each unoriented geodesic
being counted twice, once for each orientation.
So when we say, for example, that the number of geodesics of
length at most $l$ is asymptotically $\frac{e^l}{l}$,
we're talking about oriented geodesics;
the asymptotic number of unoriented geodesics would be $\frac{e^l}{2l}$.

According to the Selberg trace formula, almost conjugate hyperbolic manifolds
are \emph{isospectral}:  They have the same Laplace eigenvalues with the
same multiplicity.
(Cf.\ Randol's chapter in \cite{chavel:eigen};
Gangolli \cite{gangolli:selberg}; B{\'e}rard-Bergery \cite{bb:selberg}.)
Such manifolds can be constructed using
the well-known method of 
Sunada
\cite{sunada:method}.
Sunada's method is very flexible and powerful,
and works in contexts that go well beyond the kind of
locally homogeneous and isotropic case we encounter here
in discussing hyperbolic manifolds.
It readily yields examples of non-isometric
hyperbolic manifolds of dimension 2 and higher
that are almost conjugate and hence isospectral:
For an exposition, see Buser \cite{buser:book}.
The examples Buser describes are all orientable,
but with trivial modifications the constructions can be made to
yield isospectral pairs of non-orientable surfaces.
Since we don't know of any handy reference for this,
we'll elaborate on this in Section \ref{sunada}.

We've said that according to the Selberg trace formula,
almost conjugate manifolds are isospectral.
All the examples coming from Sunada's method
are automatically almost conjugate.
This raises the obvious question:

\begin{question} \label{question1}
If two hyperbolic manifolds are isospectral,
must they be almost conjugate?
\end{question}

In the case of orientable surfaces,
where there is no twisting to contend with,
the answer is yes:
This is `Huber's Theorem'
\cite{huber:selbergI},
dating back to 1959.
Nowadays we recognize this as a
direct consequence of the Selberg trace formula,
which shows how to `read off' from the spectrum
the lengths of the geodesics
(and, it should go without saying, the associated multiplicities).

The purpose of this paper is to prove that the answer is
still yes for surfaces, even without the orientability assumption.

\begin{thm} \label{thm1}
If two hyperbolic surfaces (not necessarily orientable or
even connected) are isospectral, then they are almost conjugate.
\end{thm}

Here, in contrast to the orientable case,
we cannot use the Selberg trace formula
to read off the lengths of geodesics directly from the spectrum.
In fact, as we will see, there are `scenarios' for constructing
counter-examples consistent with the Selberg trace formula.
But the Prime Geodesic Theorem comes to our rescue,
because we can show that any scenario for constructing a counter-example
requires
the frequent participation of a large number of geodesics of length exactly
$l$ (specifically, at least a constant times $\frac{e^l}{l}$),
and having this many geodesics of the same length is forbidden by
the Prime Geodesic Theorem.

In higher dimensions Question \ref{question1} remains open, even in the 
case of connected orientable manifolds.  The issues at stake in higher
dimensions are well illustrated in the 
proof of Theorem \ref{thm1}---so you might be interested in this theorem even if
you don't see why anyone would care about non-orientable surfaces.

To see that the possible
existence of isospectral hyperbolic manifolds that
are not almost conjugate is a question that must be taken seriously,
we note that among \emph{flat} manifolds,
in dimension $d \geq 3$
there exist isospectral pairs that are not almost conjugate.
Our favorite example of this is the 3-manifold pair
`Tetra and Didi' \cite{doylerossetti:cosmic}.
In the flat case,
some care is needed in defining almost conjugacy,
because while in a hyperbolic manifold geodesics come only in isolation,
in a flat manifold geodesics come in parallel families of varying
dimension.  So in the flat case, matching geodesics between
manifolds involves measuring, rather than just counting.
But Tetra and Didi will fail to be almost conjugate according to
any definition.

{\bf Note.}
For further insight into the possible existence of isospectral spaces
that are not almost conjugate, it is natural to expand the class of
spaces we're considering from manifolds to orbifolds.
(Cf.\ Dryden \cite{drydenfinite};
Dryden and Strohmaier \cite{drydenstrohmaier};
Dryden, Greenwald, and Gordon \cite{drydenetal})
Of course we need to extend the definition of `almost conjugacy'
appropriately.
We don't propose to discuss orbifolds in detail here, but 
for the benefit of those familiar with orbifolds,
we have appended some comments in Section \ref{orb} below.
Briefly, what we find is this:
Theorem \ref{thm1} extends to rule out examples among hyperbolic 2-orbifolds.
However, there are examples of isospectral
flat 2-orbifolds (necessarily disconnected) that are not almost conjugate.
And we still don't know what happens in the hyperbolic case in
dimension $\geq 3$.

\section{Outline}

Let $M$ be a hyperbolic surface.
Given a geodesic $\gamma$ of $M$,
let $l(\gamma)$ denote its length and
$\nu(\gamma)$ its imprimitivity index
(the number of times $\gamma$ runs around a primitive ancestor).
Define the \emph{weight} $\wt(\gamma)$ as follows:
\begin{equation}
\wt_M(\gamma)=
\dilemma{
\frac{1}{\nu} &
\mbox{if $\gamma$ is orientation-preserving}; \\
\\
\frac{1}{\nu} \tanh(l/2) &
\mbox{if $\gamma$ is orientation-reversing}.
}
\end{equation}
Define the \emph{total weight function}
\begin{equation}
\tw_M(l) = \sum_{l(\gamma)=l} \wt(\gamma)
.
\end{equation}

From the work of B\'{e}rard-Bergery
\cite{bb:selberg}
and Gangolli
\cite{gangolli:selberg}
we extract the following:

\begin{prop} \label{prop1}
Let $M$ and $N$ be hyperbolic surfaces, possibly non-orientable
or disconnected.
$M$ and $N$ are isospectral if and only if
$\tw_M = \tw_N$.
\end{prop}

We outline the proof in section \ref{prop1proof}.

In light of Proposition \ref{prop1},
to prove Theorem \ref{thm1} above it will suffice to prove

\begin{thm} \label{thm2}
If $M$ and $N$ are hyperbolic surfaces and $\tw_M = \tw_N$, then
$M$ and $N$ are almost conjugate.
\end{thm}

Observe that this is a purely geometrical statement:
All reference to the Laplace
spectrum has been laundered through the total weight function.

To prove Theorem \ref{thm2},
we will analyze how we might engineer agreement between
$\tw_M$ and $\tw_N$ without having total agreement between the geodesics
of $M$ and $N$, and show that this is not possible without
having infinitely many lengths $l$ for which the number of geodesics
of length exactly $l$  is at least $C \frac{e^l}{l}$, for $C>0$.
This will contradict the following Proposition, which is a simple
consequence of the so-called `Prime Geodesic Theorem'.

\begin{prop} \label{prop2}
For any compact hyperbolic surface, the number of geodesics of
length exactly $l$ is $o(\frac{e^l}{l})$.
\end{prop}

{\bf Proof.}
According to the Prime Geodesic Theorem
(see \cite{phillipssarnak:PGT}),
for a connected hyperbolic surface
(whether orientable or not)
the number $F(l)$ of geodesics of length at most $l$ is asymptotic to
$\frac{e^l}{l}$.
The number $f(l)$ of geodesics of length exactly $l$ is given by the
jump of $F$ at $l$:
\begin{equation}
f(l) = \lim_{s \goesto l+} F(s) - \lim_{s \goesto l-} F(s)
.
\end{equation}
But if $F$ is any positive increasing function asymptotic to $G$,
the jumps of $F$ are $o(G)$.
So $f(l) = o(\frac{e^l}{l})$.
This establishes our claim for connected surfaces.
The extension to the general case is immediate, because the
$o(\frac{e^l}{l})$ estimate
holds separately on each of the finitely many connected
components.
\qed

\section{Proof of Proposition \ref{prop1}} \label{prop1proof}

Proposition \ref{prop1}
is the 2-dimensional case of a general result that applies to
hyperbolic manifolds of any dimension
(cf. B\'{e}rard-Bergery \cite{bb:selberg})
and indeed to quotients of any rank-1 symmetric space
(cf. Gangolli 
\cite{gangolli:selberg}).
All these results grow from the original work of
Huber
\cite{huber:selbergI},
who proved Proposition \ref{prop1}
in the case of connected orientable hyperbolic surfaces,
using what amounts to a special case of the Selberg trace formula.
Huber's proof applies without essential change here, and in the more general
cases alluded to above.
Huber's paper is really beautiful and well worth reading.
For convenience we outline here how the proof goes.

In the case Huber was considering, namely orientable surfaces,
what we're calling here
the `total weight function' reduces
to what Huber calls the `length spectrum' (\emph{L\"{a}ngenspektrum}).
(Please note that this is not what various authors nowadays call the
`length spectrum'.)
Huber's Satz 7 says orientable surfaces with the same total weight function
have the same Laplace spectrum and the same volume.
(Actually he says they have the same genus, which in his case is equivalent.)
His Satz 8 says surfaces with the same Laplace spectrum have
the same weight function and the same volume.

Huber proves these results with the aid of the following
`Dirichlet series':
\[
D_M(s) =
\sum_\gamma \frac{l(\gamma)}{\nu(\gamma)}
\left ( \frac{\cosh l(\gamma)}{\cosh l(\gamma)-1} \right ) ^{\frac{1}{2}}
\cosh^{-s} l(\gamma)
.
\]
In terms of the total weight function $\tw_M$, this is
\[
D_M(s)=
\sum_l
\tw_M(l)
l
\left ( \frac{\cosh l}{\cosh l - 1} \right )^\half
\cosh^{-s} l
.
\]
Huber's entire argument goes through unscathed here,
once we have suitably modified this Dirichlet series to take proper
account of the reduced spectral contribution of orientation-reversing
geodesics,
and revised the definition of the total weight function
in a corresponding way.
The appropriate series was worked out by
B\'{e}rard-Bergery
\cite{bb:selberg}
for a hyperbolic manifold of any dimension $d$:
\[
D_M(s)=
\sum_\gamma \frac{l(\gamma)}{\nu(\gamma)}
Q(l(\gamma),A(\gamma))
\cosh^{-s} l(\gamma)
,
\]
where $A(\gamma) \in O(d-1)$ represents the action on the normal bundle of
translation around $\gamma$,
and
\[
Q(l,A)=
\left | \det \left ( I-\frac{1}{\cosh l}
\frac{A+A^\trans}{2} \right ) \right | ^{-\frac{d-1}{2}}
.
\]
This Dirichlet series converges when $s= \sigma + i t, \sigma>d-1$.
Later Gangolli
\cite{gangolli:selberg}
generalized
this even further,
to handle the case of any rank-1 symmetric space.

Specializing to surfaces,
for an orientation-preserving geodesic ($A=[1]$)
the $Q$ factor is
\[
Q(l,[1])=
\left ( 1-\frac{1}{\cosh l} \right )^{-\half}
=
\left ( \frac{\cosh l}{\cosh l - 1} \right )^\half
,
\]
while for an orientation-reversing geodesic ($A=[-1]$)
it is
\[
Q(l,[-1])=
\left ( 1+\frac{1}{\cosh l} \right )^{-\half}
=
\left ( \frac{\cosh l}{\cosh l + 1} \right )^\half
.
\]
The ratio of the latter to the former is
\[
\left ( \frac{\cosh l -1}{\cosh l + 1} \right )^\half
=
\tanh \frac{l}{2}
,
\]
which is just the factor incorporated into our definition of
the weight.
With this definition,
the Dirichlet series becomes
\[
D_M(s) =
\sum_\gamma
\wt(\gamma)
l(\gamma)
\left ( \frac{\cosh l(\gamma)}{\cosh l(\gamma)-1} \right ) ^{\frac{1}{2}}
\cosh^{-s} l(\gamma)
.
\]
In terms of the total weight function, we once more get
\[
D_M(s)=
\sum_l
\tw_M(l)
l
\left ( \frac{\cosh l}{\cosh l - 1} \right )^\half
\cosh^{-s} l
.
\]

Huber's argument, which goes through here just as in the orientable case,
is that this Dirichlet series is spectrally determined, and when restricted
to a vertical line $s= \sigma + i t$, $\sigma>1$ it represents an almost
periodic function (superposition of a discrete set of sinusoids),
whose Fourier coefficients tell the total weight function
$\tw_M$.  So the spectrum determines the total weight (and the volume).
The converse goes by observing that
$D_M(s)$ has a meromorphic extension, from the poles of which we can read off
the Laplace spectrum of $M$, together with the volume.  So the total weight
determines the spectrum.
$\qed$

\section{Proof of Theorem \ref{thm2}}

Let $\alpha_M(l)$ denote the number of primitive orientation-preserving
geodesics in $M$ of length exactly $l$, and $\beta_M(l)$ the number
of primitive orientation-reversing geodesics.

Fix two surfaces $M$ and $N$ with $\tw_M = \tw_N$, and set
\begin{equation}
a(l) = \alpha_M(l)-\alpha_N(l)
;
\end{equation}
\begin{equation}
b(l) = \beta_N(l) - \beta_M(l)
.
\end{equation}
Note that, in the second definition, $M$ and $N$ have traded places.

Our job is to show that $a(l)=b(l)=0$ for all $l$.
The condition $\tw_M = \tw_N$ tells us that
\[
\sum_{k \in \Z^+} \frac{1}{k} a(\frac{l}{k}) =
\sum_{\mbox{$k$ odd}} \frac{1}{k} b(\frac{l}{k}) \tanh(l/(2k)) +
\sum_{\mbox{$k$ even}} \frac{1}{k} b(\frac{l}{k})
.
\]
Note how on the right-hand side we have had to distinguish between
odd and even $k$,
since going around an
orientation-reversing geodesic an even number of times yields
an orientation-preserving geodesic.

This system of constraints on the integer-valued
functions $a$ and $b$ has solutions which
at first blush look like they might permit the construction of a
counter-example.
To get the simplest solutions,
fix an integer $q \geq 2$,
and  set
\[
l_0 = \log q
,
\]
so that
\[
\tanh (\frac{nl_0}{2}) = \frac{q^n-1}{q^n+1}
.
\]
Let
\[
c_n = 
\frac{1}{n} \sum_{j|n} \mu(n/j) q^j
,
\]
where $\mu$ is the usual M\"{o}bius function:
\[
\mu(n) =
\dilemma{
(-1)^k & \mbox{if $n$ is a product of $k$ distinct primes;}
\\
0 & \mbox{otherwise.}
}
\]
We get a solution to our equations by taking $a(l)=b(l)=0$ when $l$ is not
a multiple of $l_0$, and setting
\[
a(l_0)=q-1,
\]
\[
b(l_0)=q+1,
\]
\[
a(2l_0)=1,
\]
\[
b(2l_0) =0,
\]
\[
a(n l_0) = b(n l_0) = c_n \;\;\mbox{, $n=3,5,7,\ldots$}
,
\]
\[
a(n l_0) = b(n l_0) = 0 \;\;\mbox{, $n=4,6,8,\ldots$}
.
\]

The problem with these solutions is that they grow too fast:
For $n$ odd, $a(n l_0)$ is asymptotic to $q^n/n$.
Our proof of Theorem \ref{thm2} proceeds by showing that 
this kind of runaway growth is unavoidable.

{\bf Note.}
When $q=2$, for the sequence $c_1,c_2,\ldots$ we get
\[
2, 1, 2, 3, 6, 9, 18, 30, 56, 99, \ldots
;
\]
when $q=3$ we get
\[
3, 3, 8, 18, 48, 116, 312, 810, 2184, 5880, \ldots
.
\]
Looking these sequences up in Neil Sloane's online encyclopedia
of integer sequences reveals that $c_n$ tells
the number of primitive length-$n$ necklaces with beads of $q$ colors,
when turning the necklace over is not allowed.
Or in other words, $c_n$ tells the number of equivalence classes of
the set $\{0,\ldots,q-1\}^n$ under the action of the cyclic group
$\Z/n\Z$.
When $q$ is a prime power, we have the alternative interpretation of $c_n$
as the number of irreducible monic polynomials
of degree $n$ over the field with $q$ elements.
This seems suggestive:
We can't use these solutions in the
context of isospectral hyperbolic
surfaces,
but perhaps we could get some mileage out of them in a different
context\ldots .
We leave that question for another day,
and get back to the proof of Theorem \ref{thm2}.

Let 
\begin{equation}
L = \{l \suchthat \mbox{$a(l) \neq 0$ or $b(l) \neq 0$} \}
\end{equation}
and
\begin{equation}
L_0 =
\{ l \in L \suchthat
\mbox{$l$ is not a multiple of any other element of $L$}
\}
.
\end{equation}
According to this definition, $L$ is the set of lengths of geodesics
where $M$ and $N$ exhibit different behavior,
and $L_0$ consists of those lengths which are minimal
with respect to the partial order where $l \preceq m$ means
$m = k l$, $k \in \Nonetwo$.
Every element of $L$ sits above some minimal element, i.e.
\begin{equation}
L \subseteq L_0 \Nonetwo
.
\end{equation}
Our job is to show that $L_0 = \emptyset$.

\begin{lemma} \label{lemma1}
$|L_0| < \infty$
\end{lemma}

{\bf Proof.}
Suppose $l \in L_0$.
By assumption, $\tw_M(l) = \tw_N(l)$.
Because $l$ is minimal in $L$,
any contributions by imprimitive geodesics to
$\tw_M(l)$ are exactly matched by contributions to $\tw_N(l)$.
This means that the contributions of primitive geodesics of length $l$
must match:
\begin{equation}
a(l) = \tanh(l/2) b(l)
.
\end{equation}
Assume for convenience that $a(l)>0$, and hence $b(l)>a(l)$.
Rewrite the equation above:
\begin{equation}
b(l) - a(l) = b(l) (1-\tanh(l/2))
;
\end{equation}
\begin{equation}
b(l) =
\frac{b(l)-a(l)}{1-\tanh(l/2)}
.
\end{equation}
When $l$ is large,
\begin{equation}
b(l) =
\frac{b(l)-a(l)}{1-\tanh(l/2)}
\approx
\frac{e^l}{2} (b(l)-a(l))
\geq
\frac{e^l}{2}
.
\end{equation}
According to Proposition \ref{prop2},
the total number of geodesics of length exactly $l$
is $o(\frac{e^l}{l})$.
Here we have at least something on the order of 
$\frac{e^l}{2}$ geodesics of length $l$.
This puts an upper bound on $l$,
and thus forces $|L_0|<\infty$.
\qed

Let $\Podd$ denote the set of odd primes.

\begin{lemma} \label{lemma2}
For any $l \in L_0$,
only a finite number of the odd prime multiples of $l$
are also multiples of an element of $L_0$ differing from $l$.
Specifically,
\begin{equation}
|l \Podd \intersect (L_0-\{l\})\Nonetwo| \leq |L_0|-1
.
\end{equation}
\end{lemma}

{\bf Proof.}
If $l_1 \in L_0$, $l_1 \neq l$, then
$| l \Podd \intersect l_1 \Nonetwo| \leq 1$.
(This is a simple fact about divisibility:  It has nothing
special to do with lengths of geodesics!)
\qed

Now fix any $l \in L_0$, and let $p$ be an odd prime that avoids the
finite set for which $pl \in (L_0-\{l\}) \Nonetwo$.
Since $l$ is minimal in $L$, as above we have
\begin{equation}
a(l) = \tanh(l/2) b(l)
.
\end{equation}
As above, assume for convenience that $a(l)>0$, and hence $b(l)>a(l)$.

By assumption, $\tw_M(l) = \tw_N(l)$.
The only geodesics that are `in play' at length $pl$ are those
of length $l$ or $pl$:  That was the whole point
of the restriction we've placed on $p$.
So
\begin{equation}
a(pl) + \frac{1}{p} a(l)
=
\tanh(pl/2)
\left( b(pl) + \frac{1}{p}b(l) \right)
.
\end{equation}
Let's rework this:
\begin{equation}
a(pl) - \tanh(pl/2) b(pl)
=
\frac{1}{p}
( \tanh(pl/2)b(l)-a(l) )
;
\end{equation}
\begin{equation}
a(pl) - b(pl)
+ b(pl) ( 1-\tanh(pl/2) )
=
\frac{1}{p} ( \tanh(pl/2)b(l)-a(l) )
;
\end{equation}
\begin{equation}
b(pl) =
\frac
{\frac{1}{p} ( \tanh(pl/2)b(l)-a(l) )  + (b(pl) - a(pl))}
{1-\tanh(pl/2)}
.
\end{equation}

Look at the numerator here.
For $p$ large, $\frac{1}{p}(\tanh(pl/2)b(l)-a(l))$ is close to
$\frac{1}{p} (b(l)-a(l))$, and $b(l)-a(l)$ is a positive integer.
And $b(pl)-a(pl)$ is always an integer: Not necessarily
a positive integer, just some integer.
As soon as $p$ is larger than $2(b(l)-a(l))$,
$\frac{1}{p} (b(l)-a(l))$ will be a positive fraction
smaller than $1/2$,
and adding an integer to it can only increase its absolute value.
This means that for $p$ large,
the smallest the numerator can be in absolute value
is something like $\frac{1}{p}(b(l)-a(l))$,
which is at least $\frac{1}{p}$.

Meanwhile, the denominator is $1-\tanh(pl/2) \approx 2 e^{-pl}$.
So $b(pl)$ is bigger than something like
$\frac{e^{pl}}{2 p}$.
This contradicts Proposition \ref{prop2}---unless $L_0$ is
empty!
So $L_0 = \emptyset$, and $M$ and $N$ are almost conjugate.
\qed

\section{Isospectral nonorientable surfaces} \label{sunada}

It is well known that there are many
examples of isospectral closed hyperbolic surfaces.
The first examples go back to Vigneras \cite{vigneras:isospectral},
who constructed arithmetic examples from quaternion algebras.
More recent constructions have used Sunada's method.
Sunada's method is very flexible, and can produce nonorientable
examples as easily as orientable examples.
But as we don't know of a reference for this,
we briefly outline the procedure here.
For necessary backgound, see Buser \cite{buser:book}.

If you take any 
pair of Sunada isospectral closed hyperbolic surfaces without
boundary, then they have a common quotient. Now just add
what Conway calls a `cross-handle' to this quotient,
i.e., take the connected sum with a Klein bottle.
Or more generally, take the
connected sum with any closed non-orientable surface.  Put the
hyperbolic metric on this new quotient, and lift everything (the
cross-handles and the metric) back up to the covers.  The
resulting surfaces are isospectral and nonorientable.

Another way of producing isospectral nonorientable pairs is to
change some of the gluings in known orientable examples where
isospectrality is proven using transplantation.
To take a specific example, consider the surfaces described by
Buser \cite{buser:book}, Chapter 11, page 304.
If you reinterpret Buser's gluing diagrams (Figures 11.5.1 and 11.5.2)
so that the identifications on the $\beta$ geodesics are by
translation, you get a non-orientable isospectral pair.  The
$\beta$ identifications now add four cross-handles, rather than
four handles.
The transplantation method proving isospectrality in the
orientable case continues to work here as well.

\section{Comments on orbifolds} \label{orb}

Here, as promised above, are some brief comments about orbifolds.

There are three independent examples of
isospectral flat (disconnected) 2-orbifolds that are not almost conjugate,
one involving quotients of a square torus, and two involving
quotients of a hexagonal torus.
We describe them using Conway's
orbifold notation \cite{conway:orbifold}.

A standard square torus has as 2- and 4-fold quotients a $2222$
orbifold and a $244$ orbifold.  If we call the torus $S_1$ and the quotients
$S_2$ and $S_4$, spectrally
\begin{equation}
S_1 + 2 S_4 = 3 S_2
,
\end{equation}
i.e., you can't hear the
difference between a torus with two $244$s, and a trio of $2222$s.

A standard hexagonal torus $H_1$ has as 2-, 3-, and 6-fold quotients
a $2222$ orbifold $H_2$ (this is a regular tetrahedron); a $333$
orbifold $H_3$; and a $236$ orbifold $H_6$.
Spectrally,
\begin{equation}
H_2 + H_6 = 2 H_3
\end{equation}
and
\begin{equation}
H_1 + H_3 + H_6 = 3 H_2
.
\end{equation}
From these relations we can derive, for example:
\begin{equation}
H_1 + 3 H_3 = 4 H_2;
\end{equation}
\begin{equation}
H_1 + 4 H_6 = 5 H_3;
\end{equation}
\begin{equation}
2 H_1 + 3 H_6 = 5 H_2
.
\end{equation}

These examples arise from a careful analysis of the contributions
of rotations of various orders to the spectrum via the
Selberg trace formula.
To explain just how this works would take us too far afield.
However, it is possible to verify isospectrality in these examples
in a direct and elementary way by using Fourier series to represent
explicitly the eigenfunctions of the component orbifolds, and checking
that eigenvalues match up.

Among hyperbolic 2-orbifolds,
no such examples exist,
whether connected or not.
This is a corollary of Theorem \ref{thm1}, together with the observation that,
in contrast to the flat case,
in the hyperbolic case elliptic elements of differing order make
distinguishable contributions to the spectrum.
Again, to go further into detail would take us too far afield.

\section*{Thanks}

We're grateful to Emily Proctor, Carolyn Gordon, and Emily
Dryden for stimulating discussions,
and to Humboldt University
(specifically the Sonderforschungsbereich 647 project)
for hospitality and support during the writing of this paper.
We're sorry to have to thank Peter Sarnak for pointing out
that if a function is approaching a limit,
its jumps must be getting smaller and smaller.

\end{document}